\documentclass[10pt]{article}
\usepackage[mathscr]{eucal}
\usepackage{algpseudocode}

\usepackage{algorithm}
\usepackage{color}
\usepackage{amsmath}
\usepackage{amssymb}
\usepackage{graphicx}

\def\R{{\mathbb R}}

\def\N{{\mathbb N}}

\def\Z{{\mathbb Z}}

\def\OO{\mathbb{O}}

\def\HOT{\mathrm{HOT}}

\def\s{\sigma}
\def\f{\frac}
\def\p{\partial}
\def\q{\quad}

\def\la{\langle}
\def\ra{\rangle}

\def\a{{\boldsymbol a}}

\def\bb{{\boldsymbol b}}

\def\bsd{{\boldsymbol d}}

\def\bsf{{\boldsymbol f}}
\def\bff{{\mathbf f}}

\def\bi{{\mathbf i}}
\def\bsi{{\boldsymbol i}}

\def\bk{{\boldsymbol k}}

\def\bu{\boldsymbol{u}}
\def\bsv{\boldsymbol{v}}

\def\x{\boldsymbol{x}}

\def\A{{\mathbf A}}

\def\I{{\mathbf I}}

\def\bsP{{\boldsymbol P}}
\def\bR{{\boldsymbol R}}

\def\W{{\mathbf W}}

\def\bsT{\boldsymbol{T}}

\def\mA{{\mathcal A}}

\def\mH{{\mathcal H}}

\def\mN{{\mathcal N}}

\def\mT{{\mathcal T}}

\def\msX{{\mathscr X}}

\def\eep{\boldsymbol{\eps}}

\def\bi{\begin{itemize}} \def\ei{\end{itemize}}
\def\be{\begin{eqnarray*}}
\def\ee{\end{eqnarray*}}

\def\0{{\mathbf 0}}

\def\pv{\mathrm{pv}}

\newcommand{\beq}{\begin{equation}}
\newcommand{\eeq}{\end{equation}}

\def\bth{\boldsymbol{\theta}}
\def\La{{\mathbf\Lambda}}

\def\Ps{\mathbf{\Psi}}
\def\OOm{\mathbf{\Omega}}
\def\wt{\widetilde}

\def\eref#1{(\ref{#1})}

\newcommand{\eps}{\varepsilon}
\def\la{\langle}
\def\ra{\rangle}

\newcommand{\xx}{\mathbf{x}}

\newcommand{\dis}{\displaystyle}

\def\XXint#1#2#3{{\setbox0=\hbox{$#1{#2#3}{\int}$ }
\vcenter{\hbox{$#2#3$ }}\kern-.55\wd0}}

\begin{document}

\markboth{J. K. Choi AND B. Dong AND X. Zhang}
         {Tight Frame based Limited Tomography}

\title{Limited Tomography  Reconstruction via Tight Frame and  Sinogram Extrapolation  }

\author{Jae Kyu Choi\thanks{Institute of Natural Sciences, Shanghai Jiao Tong University. Email: jaycjk@sjtu.edu.cn} \qquad
Bin Dong \thanks{Beijing International Center for Mathematical Research Peking University. Email:dongbin@math.pku.edu.cn} \qquad
Xiaoqun Zhang\thanks{Institute of Natural Sciences, School of Mathematical Sciences, and MOE-LSC, Shanghai Jiao Tong University. Email: xqzhang@sjtu.edu.cn}}



\maketitle

\begin{abstract}
X-ray computed tomography (CT) is one of widely used diagnostic tools for medical and dental tomographic imaging of the human body.  However, the standard filtered backprojection reconstruction method requires the complete knowledge of the projection data.  In the case of limited data, the inverse problem of CT becomes more ill-posed, which makes the reconstructed image deteriorated by the artifacts.  In this paper, we consider two dimensional CT reconstruction using the horizontally truncated projections.  Over the decades, the numerous results including the sparsity model based approach has enabled the reconstruction of the image inside the region of interest (ROI) from the limited knowledge of the data.  However, unlike these existing methods, we try to reconstruct the entire CT image from the limited knowledge of the sinogram via the tight frame regularization and the simultaneous sinogram extrapolation.  Our proposed model shows more promising numerical simulation results compared with the existing sparsity model based approach.
\end{abstract}
\textbf{Classification}: 65N20, 65N21, 94A08 

\medskip
\section{Introduction}

X-ray computed tomography (CT) is a widely used diagnostic tool for medical and dental tomographic imaging of the human body.  It provides tomographic images of the human body by assigning an X-ray attenuation coefficient to each pixel  \cite{S.Tohnak2007}.  Let $u$ denote the (unknown) image to be reconstructed.  We further assume that $u$ is supported in the unit ball $B(\0,1)$ in $\R^2$.  In the case of two dimensional parallel beam CT, the projection data (or sinogram) $f$ for each $\varphi\in[0,2\pi)$ and $s\in\R$ is obtained via the following Radon transform \cite{Radon1917}:
\begin{align}\label{Radon_Transform}
f(\varphi,s)=P u(\varphi,s)=\int_{-\infty}^{\infty}u(s\bth+t\bth^{\perp})dt
\end{align}
where $\bth=(\cos\varphi,\sin\varphi)$ and $\bth^{\perp}=(-\sin\varphi,\cos\varphi)$.   Then the tomographic image $u$ is reconstructed via the following Radon inversion formula \cite{R.N.Bracewell1967,Natterer2001}:
\begin{align}\label{FBP}
u(\x)=\f{1}{2\pi^2}\int_0^{\pi}\int_{-\infty}^{\infty}\f{1}{\x\cdot\bth-s}\left[\f{\p}{\p s}Pu(\varphi,s)\right]dsd\varphi
\end{align}
with $\x=(x_1,x_2)\in\R^2$.

\begin{figure}[h]
\begin{center}
\begin{tabular}{cccc}
\begin{minipage}{3cm}
\includegraphics[width=3cm]{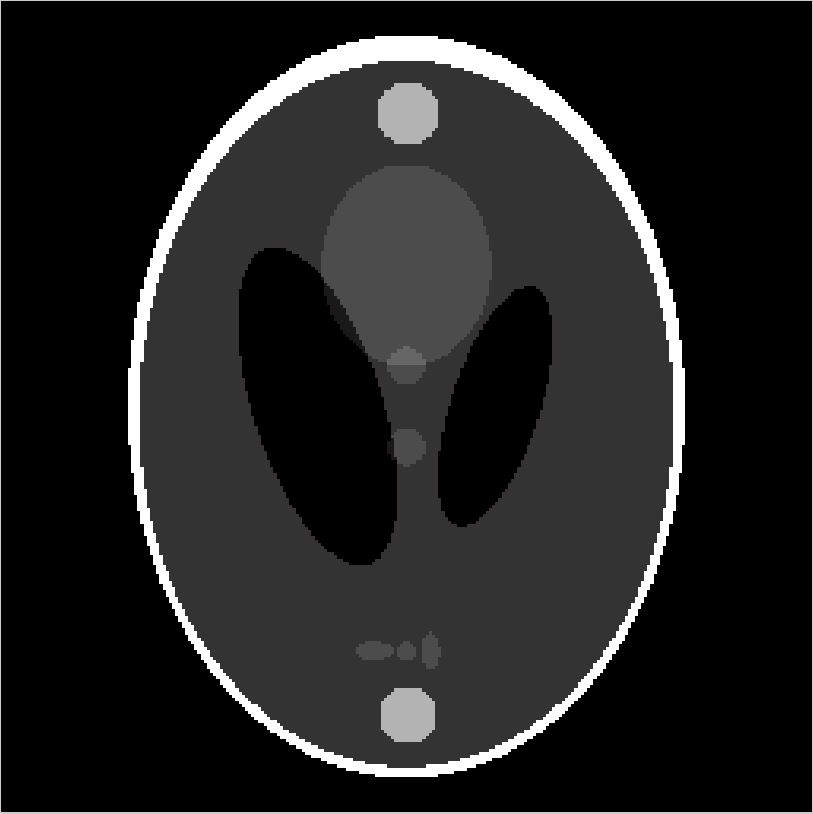}
\end{minipage}&
\begin{minipage}{3cm}
\includegraphics[width=3cm]{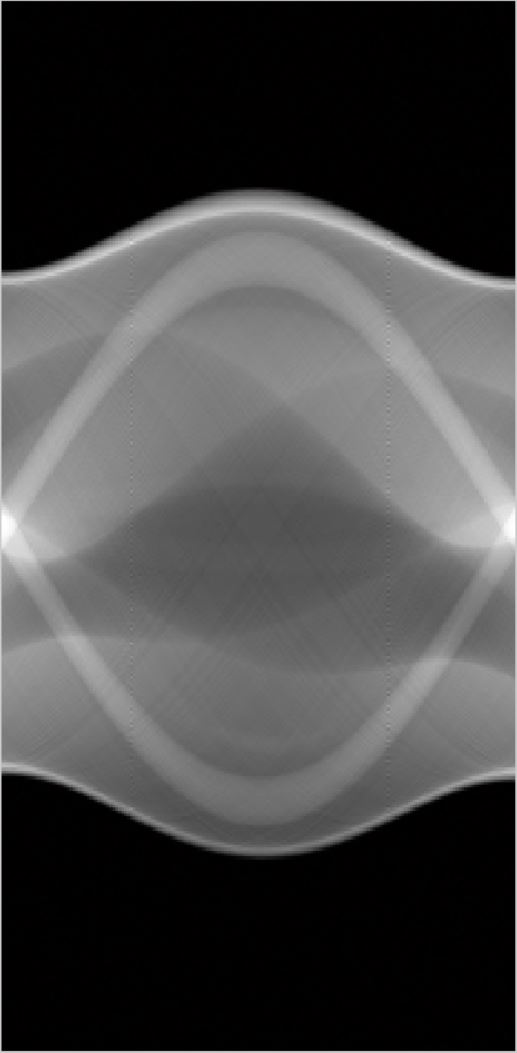}
\end{minipage}&
\begin{minipage}{3cm}
\includegraphics[width=3cm]{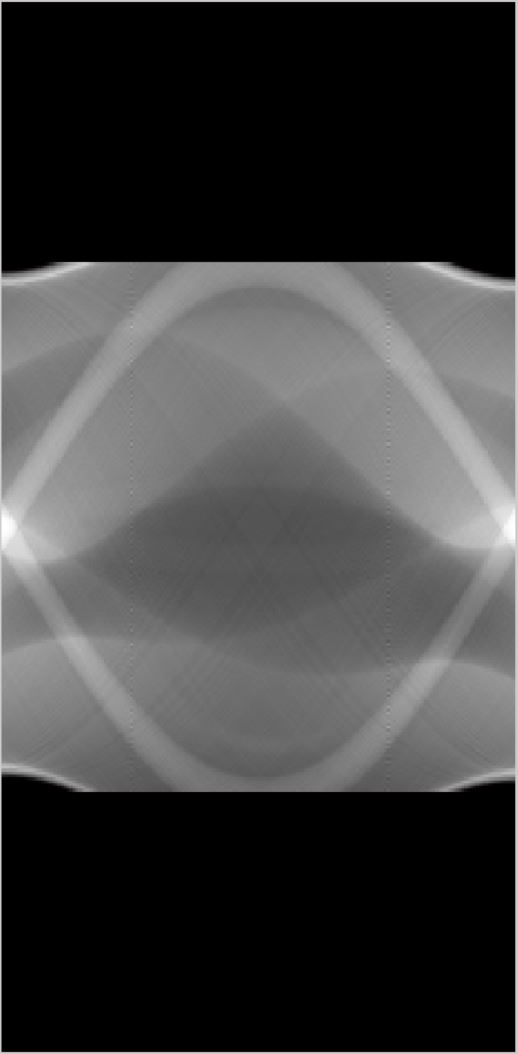}
\end{minipage}&
\begin{minipage}{3cm}
\includegraphics[width=3cm]{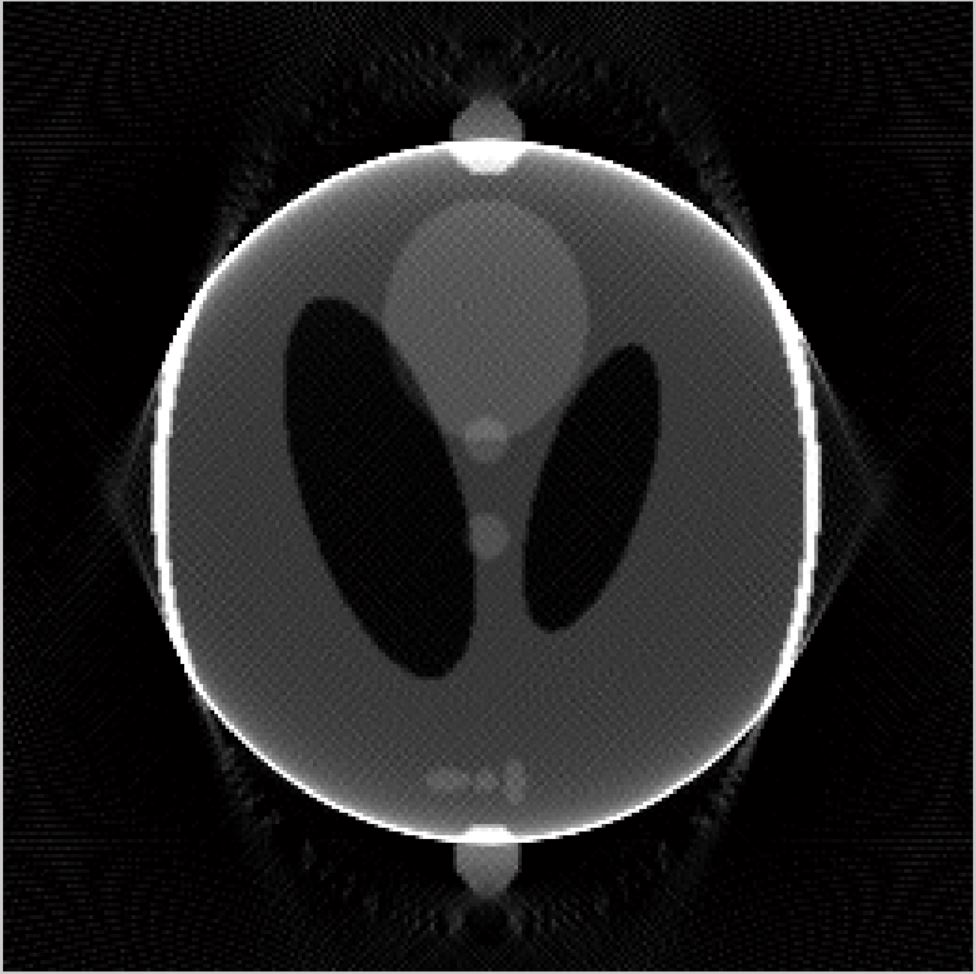}
\end{minipage}\\
(a) & (b) & (c) & (d)
\end{tabular}
\caption{The challenging issues in the limited tomography which is illustrated by (a) the original phantom image, (b) the full projection data covering $[0,\pi)\times\R$, (c) the limited data due to the detector size, and (d) the reconstructed image by \eref{FBP} with the limited data.  As we can see, the reconstructed image is corrupted by the ambiguity of $P$.}\label{issue_data_truncation}
\end{center}
\end{figure}

However, the limitation of \eref{FBP} lies in the fact that it requires so-called \emph{complete data} \cite{Klann2011,Quinto1993,Quinto2006}, which means that the measured projection data should cover at least the range $(\varphi,s)\in[0,\pi)\times\R$.  In the case of \emph{limited data} where the sinogram $f(\varphi,s)$ is measured only on a subset of $[0,\pi)\times\R$ due to the reduced size of detector \cite{M.Lee2015,G.Wang2013,J.P.Ward2015} and/or the reduced number of projections \cite{K.Choi2010,E.Y.Sidky2008,X.Pan2009}, the reconstruction is more ill-posed than the complete data case \cite{Finch1985,Klann2011}.  In particular, if the projection data is available only on $\Lambda_{\mu}\subseteq[0,\pi)\times\R$ given as
\begin{align}
\Lambda_{\mu}=\{(\varphi,s)\in[0,\pi)\times\R:|s|<\mu<1\},
\end{align}
then there exists a nontrivial function $g$ called the amgibuity of $P$ \cite{G.Wang2013}, i.e. $P g=0$ in $\Lambda_{\mu}$ \cite{G.Wang2013,J.Yang2010}.  Since it has been proven that $g$ is nonconstant in the region of interest (ROI) $B(\0,\mu)$ \cite{Natterer2001,G.Wang2013}, the reconstructed image via \eref{FBP} with the projection data $f$ restricted on $\Lambda_\mu$ will be deteriorated by this ambiguity, as shown in Fig. \ref{issue_data_truncation}.

During the past decades, the development of CT theories has enabled the unique and stable reconstruction of CT image from the limited knowledge of the sinogram \cite{M.Courdurier2008,W.Han2009,H.Kudo2008,M.Lee2015,J.P.Ward2015,J.Yang2010,Y.Ye2007,H.Yu2009}.  These existing reconstruction methods mainly focus on the reconstruction of the image inside $B(\0,\mu)$.  However, since the sinogram restricted to $\Lambda_{\mu}$ may not necessarily agree with the projection of the ROI part of image due to the integral geometry of the Radon transform \cite{L.C.B.dosSantos2014}, we may have to solve another inverse problem to remove the projections of the external contribution for the reconstruction of image inside the ROI.

In this paper, we propose a new image reconstruction method for the limited tomography.  We note that with the advent of compressed sensing theory \cite{A.M.Bruckstein2009,Cand`es2006,E.J.Cand`es2006,Donoho2006}, the sparsity model based reconstruction methods are proposed to reconstruct the image on the ROI from the knowledge of limited projection data \cite{W.Han2009,J.Yang2010,H.Yu2009}.  Unlike these existing sparsity model based methods that aim to reconstruct the image on the ROI only, our proposed method aims to reconstruct the entire CT image from the knowledge of the projection data given only on $\Lambda_{\mu}$.  Motivated by the idea of sinogram inpainting in \cite{B.Dong2013a}, our method is based on the reconstruction of CT image and the simultaneous extrapolation of the sinogram.  Since both the CT images and the proejection data are discrete in practice, we use the tight frames as the tool of limited tomography image reconstruction.  The numerical simulation shows that our proposed method has a possibility to recover the image outside ROI, which is believed to be difficult to reconstruct by FBP \eref{FBP} and existing methods.

{In \cite{M.Burger2014}, simultaneous total variation regularization on the image space and sinogram is also adopted for PET reconstruction. An interpolation method for sinogram data is also proposed in \cite{M.Kalke2014}.  More recently, a joint spatial-Radon domain CT image reconstruction model based on data-driven tight frames (SRD-DDTF) was proposed in \cite{R.Zhan2016}.  It should be noted that the simultaneous regularization on both CT image and sinogram using tight frame is proposed in \cite{B.Dong2013a, R.Zhan2016} in the limited angle tomography problem.   However, the use of this simultaneous regularization is so far as we know the first try in the limited tomography problem to reconstruct the entire CT image.}

The rest of this paper is organized as follows: in the following section, we briefly introduce the existing methods on the limited tomography, most of which are given in the continuous setting.  In section \ref{section3}, we introduce our tight frame based reconstruction algorithm for the limited tomography in the discrete setting.  Then the numerical simulation results are presented in section \ref{section4}.

\section{Preliminaries and Existing Methods}\label{section2}

In this section, we briefly introduce existing results on the limited tomography.  For more details, the readers may refer to \cite{G.Wang2013} and the references therein.

In the past decades, numerous extensive studies have been introduced to remove the ambiguity caused by the limited knowledge of the projection data \cite{G.Wang2013}.   These studies mainly focus on the reconstruction of $u$ in the ROI $B(\0,\mu)$ from the knowledge of the projection data $f$ in the region $\Lambda_{\mu}$.  Their results can be classified into the following two classes: the known subregion based approach \cite{M.Courdurier2008,H.Kudo2008,Y.Ye2007} and the sparsity model based approach \cite{W.Han2009,M.Lee2015,J.P.Ward2015,J.Yang2010,H.Yu2009}.  The known subregion based approach is based on the following identity
\begin{align}\label{Differentiation_Backprojection}
\mH_{\bth_0}u(\x)=-\f{1}{2\pi}P_{\bth}^*\left[\f{\p}{\p s}P u\right](\x)
\end{align}
where $\bth_0=(\cos\varphi_0,\sin\varphi_0)$ for some fixed angle $\varphi_0$.  Here, $P_{\bth_0}^*$ and $\mH_{\bth}$ respectively denote the backprojection operator and the directional Hilbert transform which are defined by
\begin{align}
P_{\bth}^*f(\x)&=\int_{\varphi_0-\pi/2}^{\varphi_0+\pi/2}f(\varphi,\x\cdot\bth)d\varphi \label{backprojection}\\
\mH_{\bth}u(\x)&=\f{1}{\pi}\pv\int_{-\infty}^{\infty}\f{u(\x-s\bth_0)}{s}ds \label{directional_Hilbert}
\end{align}
with $\pv$ denoting the Cauchy principal value \cite{E.M.Stein2011}.

From the above identities, the problem is reduced to reconstruct $u$ in the ROI $B(\0,\mu)$ from the knowledge of $\mH_{\bth_0}u$ in $B(\0,\mu)$.  Since the directional Hilbert transform is equivalent to one dimensional Hilbert transforms along the straight lines \cite{J.P.Ward2015}, this reduces the problem to the reconstruction of one dimensional function $F$ on the interval $[a,b]$ from the knowledge of its Hilbert transform on this interval.  Then the known subregion based approach adopts the analyticity of the ambiguity of Hilbert transform $\mH$ on $(a,b)$ \cite{M.Courdurier2008,H.Kudo2008,Y.Ye2007}; from the given knowledge of $F$ in an open interval contained in $(a,b)$, we can remove the ambiguity by the analytic continuation.  However, due to the numerical unstability of analytic continuation, the numerical implementation mostly relies on the projection onto convex sets method \cite{M.Defrise2006}, whose computational cost is relatively expensive \cite{M.Lee2013,G.Wang2013}.

\begin{figure}[h]
\begin{center}
\includegraphics[width=9cm]{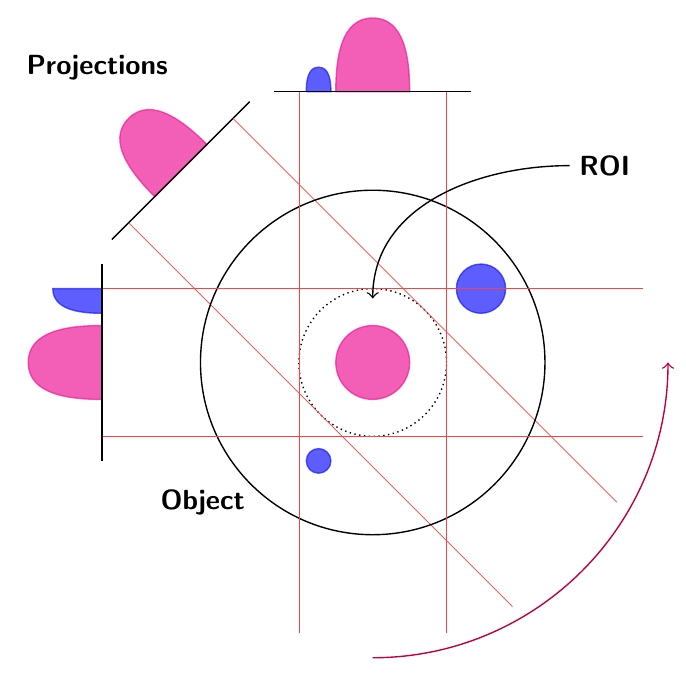}
\caption{Integral geometry of limited tomography for two dimensional parallel beam CT case.  Even though the projection data is limited due to the size of detector, the data can contain the outside ROI contribution.}\label{illustration_problem}
\end{center}
\end{figure}

Even though the unique reconstruction by known subregion based approach is guaranteed mathematically, this approach may not be practically applicable since there would be situations where no precise prior information is provided on any subregion \cite{G.Wang2013}.  Sparsity model based approach is using a regularization term that enforces the regularity of $u$: for a given limited projection $f$, the best image $u$ is found using the following constraints:
\begin{align}\label{Sparsity_Method}
\min~R(u)~~~~~\text{s.t.}~~~P u=f~~\text{in}~\Lambda_{\mu}.
\end{align}
For example, we can choose the total variation as a regularization term $R(u)$ to capture the discontinuities \cite{G.Wang2013}.  The authors in \cite{J.Yang2010} proposed a general model by using the following higher order total variation as $R(u)$:
\begin{align*}
\HOT_k(u)=\sum_{i=1}^m&\sum_{j>i,j\in\mN_i}\int_{\Gamma_{i,j}}|u_i-u_j|ds\\
&+\int_{B(\0,\mu)}\min\left\{\sqrt{\sum_{r=0}^{k+1}\left(\f{\p^{n+1}u}{\p x_1^r\p x_2^{k+1-r}}\right)^2},\sqrt{\left(\f{\p u}{\p x_1}\right)^2+\left(\f{\p u}{\p x_2}\right)^2}\right\}d\x.
\end{align*}
Here, $m$ is the number of subregions, $\Gamma_{i,j}$ is the piecewise smooth boundaries between $i$th and $j$th subregions, $u_i$ is the trace of $u$ from the $i$th subregion to $\Gamma_{i,j}$, and $\mN_i$ is the neighborhood regions of the $i$th subregion \cite{J.Yang2010}.  Under this setting, it has been shown that $u$ in $B(\0,\mu)$ can be uniquely and stably reconstructed by solving \eref{Sparsity_Method} with $R(u)=\HOT_k(u)$ provided that it is piecewise polynomial in $B(\0,\mu)$ \cite{E.Katsevich2012,J.Yang2010}.

Even though the sparsity model based approach can reconstruct the image exactly without the prior knowledge of $u$ in a subregion of $B(\0,\mu)$ \cite{J.Yang2010}, its major disadvantage is that \eref{Sparsity_Method} has to be solved for the entire image, including the exterior of $B(\0,\mu)$ \cite{L.C.B.dosSantos2014}.  One possible approach to overcome this disadvantage is to reduce the columns of discrete Radon transform matrix and consider the contribution of the tissues inside the ROI only \cite{L.C.B.dosSantos2014,K.Niinimaeki2007}.  However, it should be noted that even though the projection data $f$ is restricted on $\Lambda_{\mu}$, it contains the projections from the ROI, as well as the projections from the exterior of the ROI.  This means that the restricted projection data may not necessarily agree with the Radon transform of the ROI part of $u$, as shown in Fig. \ref{illustration_problem}.

\section{Methods and Algorithms}\label{section3}

In this section, we present our tight frame  based CT reconstruction model with a given limited projection data.  The brief explanations on the wavelet tight frame and the data driven tight frame are given in Appendix.

\subsection{Tight Frame Based CT Reconstruction Models}

Let us begin this subsection by introducing some notations.  Let $\bsf_0$ be the measured projection data defined on the $N_s\times N_p$ grid $\La$ where $N_s$ is the number of grid points on the X-ray detector and $N_p$ is the number of projections.  Note that $\La$ is the discretizaiton of $\Lambda_{\mu}$ in the Radon domain.  We denote the extrapolated projection data as $\bsf$ defined on the $M_s\times N_p$ grid $\OOm$ containing $\La$.  Let $\bR_{\La}$ denote the restriction operator defined as
\begin{align*}
\bR_{\La}\bsv[\bk]=\left\{\begin{array}{ccl}
\bsv[\bk]~&\text{if}&~\bk\in\La\\
0~&\text{if}&~\bk\notin\La.
\end{array}\right.
\end{align*}
In the noise-free case, $\bsf$ should satisfy $\bR_{\La}\bsf=\bsf_0$, which means that the extrapolated projection data should be consistent with the measured limited projection data.

Let $\bsP$ denote the discrete Radon transform matrix.  Since $\bsf$ has to be consistent with the projection $\bsP\bu$ of the reconstructed CT image $\bu$, this leads us to propose our limited tomography reconstruction model as follows:
\begin{align}\label{DTDVLT}
\min_{\bsf\geq0,0\leq\bu\leq a}~\lambda_1\|\W_1\bsf\|_1+\lambda_2\|\W_2\bu\|_1~~\text{s.t.}~~\bR_{\La}\bsf=\bsf_0,~~\bR_{\La}\bsP\bu=\bsf_0~\&~\bR_{\La^c}\bsP\bu=\bR_{\La^c}\bsf
\end{align}
with $\lambda_1$, $\lambda_2>0$.  For the noisy limited data $\bsf_0$, the first two equality constraints in \eref{DTDVLT} can be replaced with $\|\bR_{\La}\bsf-\bsf_0\|_2\leq\eps$ and $\|\bR_{\La}\bsP\bu-\bsf_0\|_2\leq\eps$ respectively, where $\eps$ is the parameter related with the standard deviation of noise.

The norm $\|\cdot\|_1$ in \eref{DTDVLT} is the isotropic $\ell^1$ norm \cite{J.F.Cai2015,B.Dong2013a} defined as
\begin{align}\label{framelet_1_norm}
\|\W\bu\|_1=\sum_{\bk\in\OO^2}\left[\sum_{l=0}^{L-1}\left(\sum_{\bsi\neq\0}\big|(\W_{l,\bsi}\bu)[\bk]\big|^2\right)^{1/2}\right].
\end{align}
with the cartesian grid $\OO^2$ where $\bu$ is defined.  Note that $\lambda_1\|\W_1\bsf\|_1$ and $\lambda_2\|\W_2\bu\|_1$ are used to guarantee the sparse representations of $\bsf$ and $\bu$ under two different tight frames $\W_1$ and $\W_2$ respectively.

In this paper, we will use two types of tight frames for comparison of the numerical simulation results.  One type that we will use is the wavelet frame system given in \cite{J.F.Cai2015,J.F.Cai2012,B.Dong2013}.  More precisely, we will use the cubic B-spline framelet system with $3$ levels of decomposition for $\W_1$ and the linear B-spline framelet system with $1$ level of decomposition for $\W_2$.  The other type of tight frame is the data driven frames adaptively learned from the preprocessed data \cite{J.F.Cai2014,B.Dong2015}. The data driven tight frame was adopted in \cite{J.Liang2014,S.Yu2015} for seismic data interpolation and CT reconstruction in \cite{W.Zhou2013}.  In our numerical simulation, the cubic B-spline framelet and the linear B-spline framelet are chosen as initial $\W_1$ and $\W_2$ respectively.

We note that \eref{DTDVLT} is more general than the analysis based approach for the limited tomography without Radon domain extrapolation:
\begin{align}\label{Analysis_LT}
\min_{0\leq\bu\leq a}~\|\W\bu\|_1~~~~~\text{s.t.}~~~~~\bR_{\La}\bsP\bu=\bsf_0.
\end{align}
When $\W$ is a B-spline framelet, then by \cite{J.F.Cai2012}, \eref{Analysis_LT} can be viewed as a finite difference approximation of sparse model based approach \eref{Sparsity_Method} with $\HOT_k(u)$ for some $k\in\N$.  In this context, we refer to \eref{Analysis_LT} with a B-spline framelet $\W$ as a sparsity model based method in the rest of this paper.  In addition, we will refer to the model \eref{DTDVLT} with B-spline framelet systems as wavelet frame based model, and \eref{DTDVLT} with data driven tight frames as data driven frame based model.

\subsection{Bregmanized Operator Splitting Algorithms}

To solve our model \eref{DTDVLT}, we note that the problem \eref{DTDVLT} is of the form
\begin{align}\label{Eq_Constraint}
\min_{\xx}~J(\xx)~~~~~\text{s.t.}~~~~\A\xx=\bff_0
\end{align}
where $\xx=[\begin{array}{c|c|c|c}
\bsf&\bu&\bsd_1&\bsd_2
\end{array}]^T$, and $J(\xx)$,  $\A$,  and ${\bff}_0$ are respectively defined as
\begin{align*}
J(\xx)=\chi_{\bsf\geq 0}(\bsf)+\chi_{0\leq\bu\leq a}(\bu)+\lambda_1\|\bsd_1\|_1+\lambda_2\|\bsd_2\|_1,\end{align*}
where $\chi (\cdot)$ is the convex set indicator function and \begin{align*}
\A=\left[\begin{array}{c|c|c|c}
\bR_{\La}&\0&\0&\0\\ \hline
\0&\bR_{\La}\bsP&\0&\0\\ \hline
\bR_{\La^c}&-\bR_{\La^c}\bsP&\0&\0\\\hline
\W_1&\0&-\I&\0\\ \hline
\0&\W_2&\0&-\I
\end{array}\right],\quad
 {\bff}_0=\left[\begin{array}{c}
\bsf_0\\ \hline
\bsf_0\\ \hline
\0\\ \hline
\0\\ \hline
\0
\end{array}\right].
\end{align*}
To solve the problem, we may use e.g. the split Bregman algorithm \cite{T.Goldstein2009,J.F.Cai2009/10} or the augmented Lagrangian method \cite{C.Wu2010}.  We use the Bregmanized operator splitting (BOS) method which is also known as split inexact  method \cite{X.Zhang2010,X.Zhang2010a} as it can transfer the constrained optimization problem into the several easy and efficient subproblems without inner iterations \cite{J.Liang2013}. The BOS scheme for solving \eref{Eq_Constraint} is as followed:
\begin{align}
\xx^{k+1}&=\arg\min_{\xx\in\mA}~\left(\mu J(\xx)+\f{\kappa}{2}\|\A\xx-\bff^{k}\|_2^2+\frac{1}{2}\|\xx-\xx^k\|_Q^2\right) \label{xupdate}\\
\bff^{k+1}&=\bff^k+\bff_0-\A\xx^{k+1}. \label{bupdate}\end{align}
Here, $\|\cdot\|_Q^2=\langle Q\cdot,\cdot\rangle$ is a semi-norm with respect to a semi-positive definite matrix $Q$ and $\bff=[\begin{array}{c|c|c|c|c}
\bsf_1&\bsf_2&\bsf_3&\bb_1&\bb_2
\end{array}]^T$.   By choosing $Q$ properly, we obtain  the following Algorithm \ref{Alg31}: for positive parameters $\kappa, \beta$
\begin{algorithm}[BOS Algorithm for \eref{DTDVLT}]\label{Alg31}\q\q
\begin{algorithmic}
\State{\textbf{Initialize:} $\bsf^{(0)}$, $\bu^{(0)}$, $\bsv_1^{(0)}$, $\bsv_2^{(0)}$, $\bsf_1^{(0)}$, $\bsf_2^{(0)}$, $\bsf_3^{(0)}$, $\bsd_1^{(0)}$, $\bsd_2^{(0)}$, $\bb_1^{(0)}$, $\bb_2^{(0)}$}
\For {$k=0$, $1$, $2$, $\cdots$}
\State{(1)  $~~\bsv_1^{k+1}=\bsf^k-\kappa^{-1}[\bR_{\La}^T(\bR_{\La}\bsf^k-\bsf_1^k)+\bR_{\La^c}^T\{\bR_{\La^c}(\bsf^k-\bsP\bu^k)-\bsf_3^k\}]$}
\State{(2)  $~~\bsv_2^{k+1}=\bu^k-\kappa^{-1}\bsP^T[\bsP\bu^k-\bR_{\La}^T\bsf_2^k-\bR_{\La^c}^T(\bsf^k-\bsf_3^k)]$}
\State{(3)  $~~\bsf^{k+1}=\arg\min_{\bsf\geq0}~\f{\beta}{2}\|\W_1\bsf-\bsd_1^k+\bb_1^k\|_2^2+\f{\kappa}{2}\|\bsf-\bsv_1^{k+1}\|_2^2$}
\State{(4)  $~~\bu^{k+1}=\arg\min_{0\leq\bu\leq a}~\f{\beta}{2}\|\W_2\bu-\bsd_2^k+\bb_2^k\|_2^2+\f{\kappa}{2}\|\bu-\bsv_2^{k+1}\|_2^2$}
\State{(5)  $~~\bsd_1^{k+1}=\arg\min_{\bsd_1}~\lambda_1\|\bsd_1\|_1+\f{\beta}{2}\|\bsd_1-\W_1\bsf^{k+1}-\bb_1^k\|_2^2$}
\State{(6)  $~~\bsd_2^{k+1}=\arg\min_{\bsd_2}~\lambda_2\|\bsd_2\|_1+\f{\beta}{2}\|\bsd_2-\W_2\bu^{k+1}-\bb_2^k\|_2^2$}
\State{(7)  $~~\bsf_1^{k+1}=\bsf_1^k+\bsf_0-\bR_{\La}\bsf^{k+1}$}
\State{(8) $~\bsf_2^{k+1}=\bsf_2^k+\bsf_0-\bR_{\La^c}\bsP\bu^{k+1}$}
\State{(9) $~\bsf_3^{k+1}=\bsf_3^k+\bR_{\La^c}(\bsP\bu^{k+1}-\bsf^{k+1})$}
\State{(10)  $~~\bb_1^{k+1}=\bb_1^k+\W_1\bsf^{k+1}-\bsd_1^{k+1}$}
\State{(11)  $~~\bb_2^{k+1}=\bb_2^k+\W_2\bu^{k+1}-\bsd_2^{k+1}$}
\EndFor
\end{algorithmic}
\end{algorithm}
The convergence condition is that $\kappa>\|\A_1^T\A_1\|$ where $\A_1=\left[\begin{array}{c|c}
\bR_{\La}&\0\\ \hline
\0&\bR_{\La}\bsP\\ \hline
\bR_{\La^c}&-\bR_{\La^c}\bsP
\end{array}\right]$.

Note that both \emph{(3)} and \emph{(4)} can be solved via the following two steps; since $\W_1$ and $\W_2$ are tight frames, we first have:
\begin{align}\label{u_closed_temp}
\bu_i^{k+1/2}&=\f{1}{\kappa+\beta_i}\big[\kappa\bsv_i^{k+1}+\beta\W_i^T(\bsd_i^k-\bb_i^k)\big]~~~~~\text{for}~~i=1,2
\end{align}
where $\bu_1=\bsf$ and $\bu_2=\bu$.  Then we have
\begin{align*}
\bsf^{k+1}&=\max(\bsf^{k+1/2},0)\\
\bu^{k+1}&=\min\big(\max(\bu^{k+1/2},0),a\big),
\end{align*}
where both $\min$ and $\max$ are componentwise operation.  In addition, \emph{(5)} and \emph{(6)} have simple analytical solutions as well:
\begin{align}\label{d_closed}
\bsd_i^{k+1}=\mT_{\lambda_i/\beta}\left(\W_i\bu_i^{k+1}+\bb_i^k\right)~~~~~\text{for}~~~i=1,2,
\end{align}
where $\mT_{\alpha}$ is the soft-thresholding operator \cite{J.F.Cai2012,B.Dong2013a} defined as
\begin{align}\label{soft_threshold}
\mT_{\alpha}(\bsv)_{l,\bsi}[\bk]=\left\{\begin{array}{lcl}
\bsv_{l,\bsi}[\bk]~&\text{if}&~\bsi=\0,
\vspace{0.4em}\\
\dis{\f{\bsv_{l,\bsi}[\bk]}{R_l[\bk]}\max\big(R_l[\bk]-\alpha,0\big)}~&\text{if}&~\bsi\neq\0,
\end{array}\right.
\end{align}
with $R_l[\bk]=\left[\sum_{\bsi\neq\0}|\bsv_{l,\bsi}[\bk]|^2\right]^{1/2}$.

For the initialization, we use $\bsf^0=\bsf_0$ and $\bu^0$ as the following least square solution
\begin{align}\label{LS_Sol}
\min_{\bu}~\|\bR_{\La}\bsP\bu-\bsf_0\|_2^2,
\end{align}
which is obtained by FBP \eref{FBP} with a measured limited data $\bsf_0$.  For the remaining variables, we initialize $\bsd_1^0=\bb_1^0=\bsd_2^0=\bb_2^0=\bsf_3^0=\0$ and $\bsf_1^0=\bsf_2^0=\bsf_0$.  For the data driven tight frame case, we first learn $\W_1$ and $\W_2$ using $\bsf$ and $\bu$ calculated from the previous iteration.  Then we perform Algorithm \ref{Alg31} with learned $\W_1$ and $\W_2$ as an inner iteration.

In addition, the sparsity model based method \eref{Analysis_LT} can be solved in a similar way using the BOS method:
\begin{align*}
\bsv^{k+1}&=\bu^k-\kappa^{-1}\bsP^T\bR_{\La}^T(\bR_{\La}\bsP\bu^k-\bsf_0^k)\\
\bu^{k+1}&=\arg\min_{0\leq\bu\leq a}~\f{\beta}{2}\|\W\bu-\bsd^k+\bb^k\|_2^2+\f{\kappa}{2}\|\bu-\bsv^{k+1}\|_2^2\\
\bsd^{k+1}&=\arg\min_{\bsd}~\mu\|\bsd\|_1+\f{\beta}{2}\|\bsd-\W\bu^{k+1}-\bb^k\|_2^2\\
\bb^{k+1}&=\bb^k+\W\bu^{k+1}-\bsd^{k+1}\\
\bsf_0^{k+1}&=\bsf_0^k+\bsf_0-\bR_{\La}\bsP\bu^{k+1}
\end{align*}
with initializations $\bsd^0=\bb^0=\0$, $\bsf_0^0=\bsf_0$, and $\bu^0$ being the solution of \eref{LS_Sol}.

\section{Numerical Simulations}\label{section4}

\begin{figure}[h]
\begin{center}
\begin{tabular}{c}
\begin{minipage}{3cm}
\includegraphics[width=3cm]{Figure/Reference_Image.JPG}
\end{minipage}
\end{tabular}\vspace{0.2em}\\
\begin{tabular}{cccc}
\begin{minipage}{3cm}
\includegraphics[width=3cm]{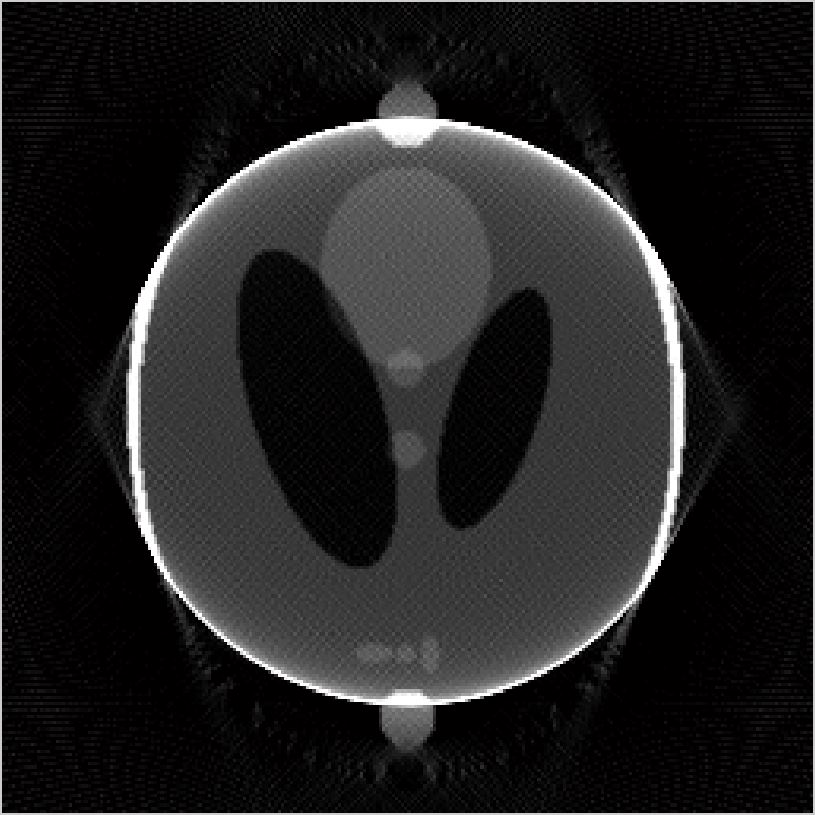}
\end{minipage}&
\begin{minipage}{3cm}
\includegraphics[width=3cm]{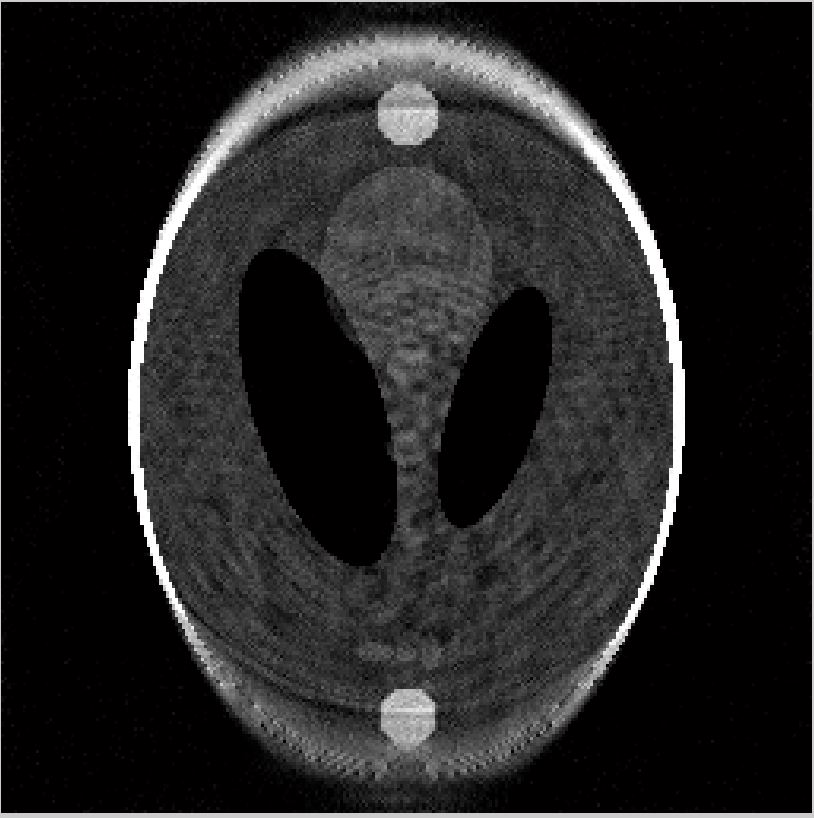}
\end{minipage}&
\begin{minipage}{3cm}
\includegraphics[width=3cm]{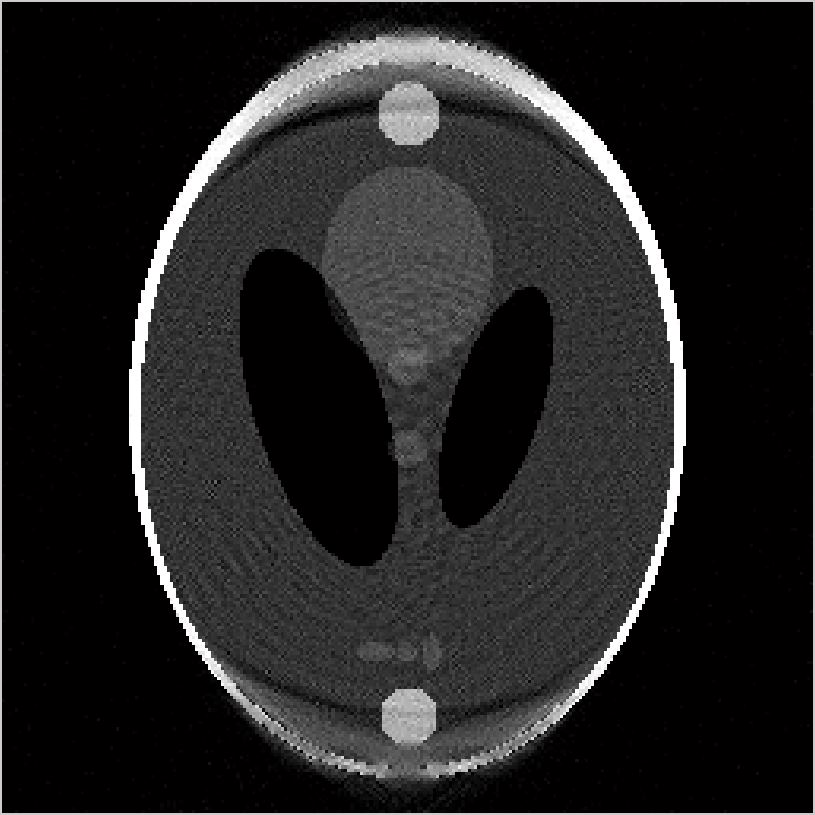}
\end{minipage}&
\begin{minipage}{3cm}
\includegraphics[width=3cm]{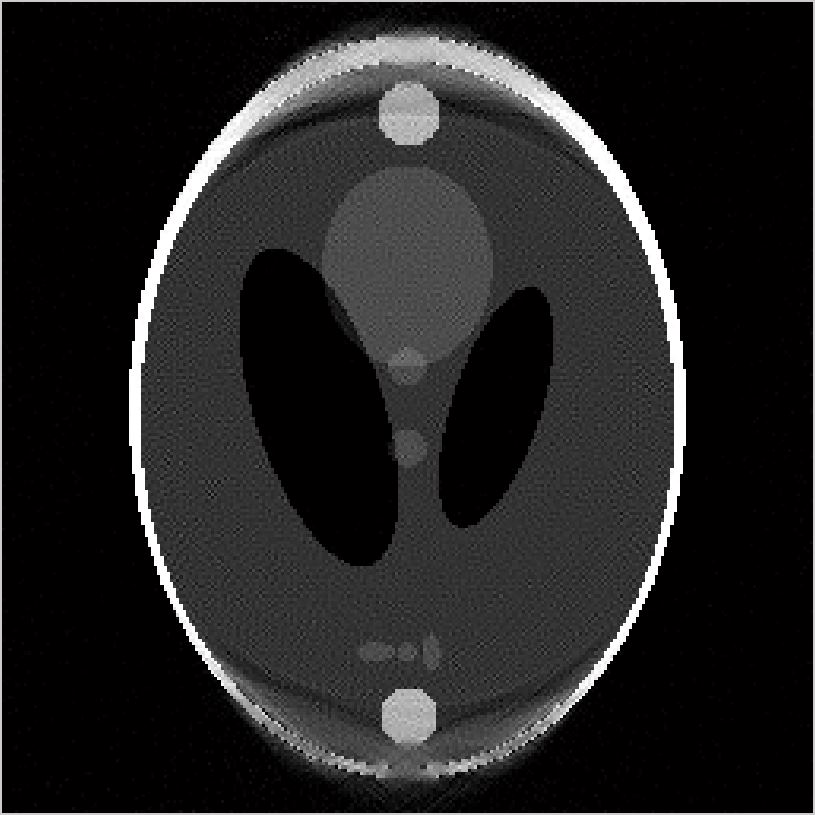}
\end{minipage}\vspace{0.3em}\\
\begin{minipage}{3cm}
\includegraphics[width=3cm]{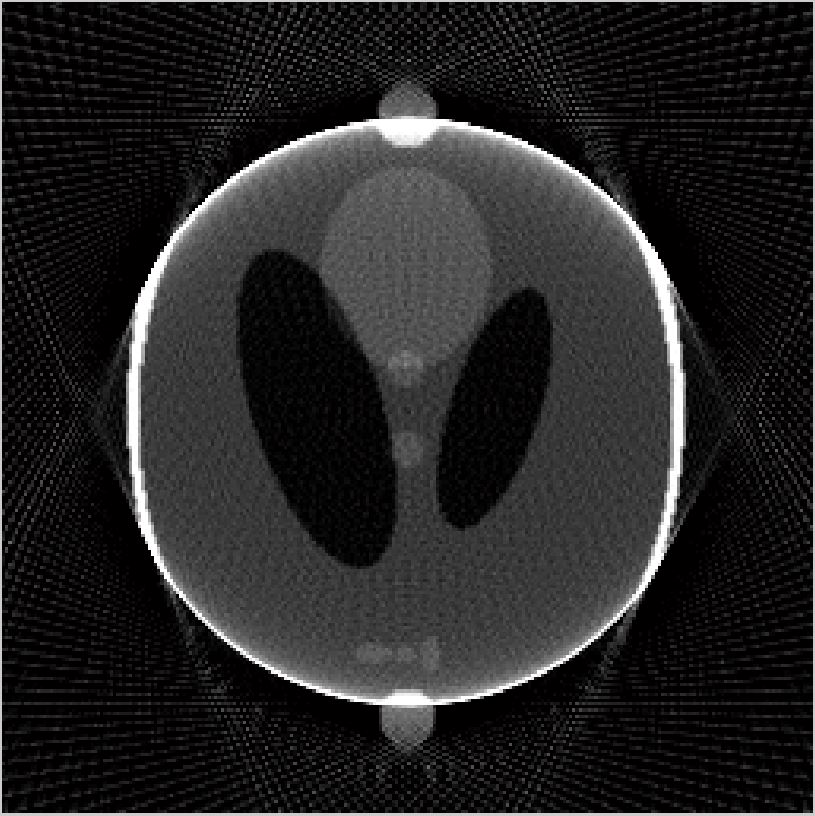}
\end{minipage}&
\begin{minipage}{3cm}
\includegraphics[width=3cm]{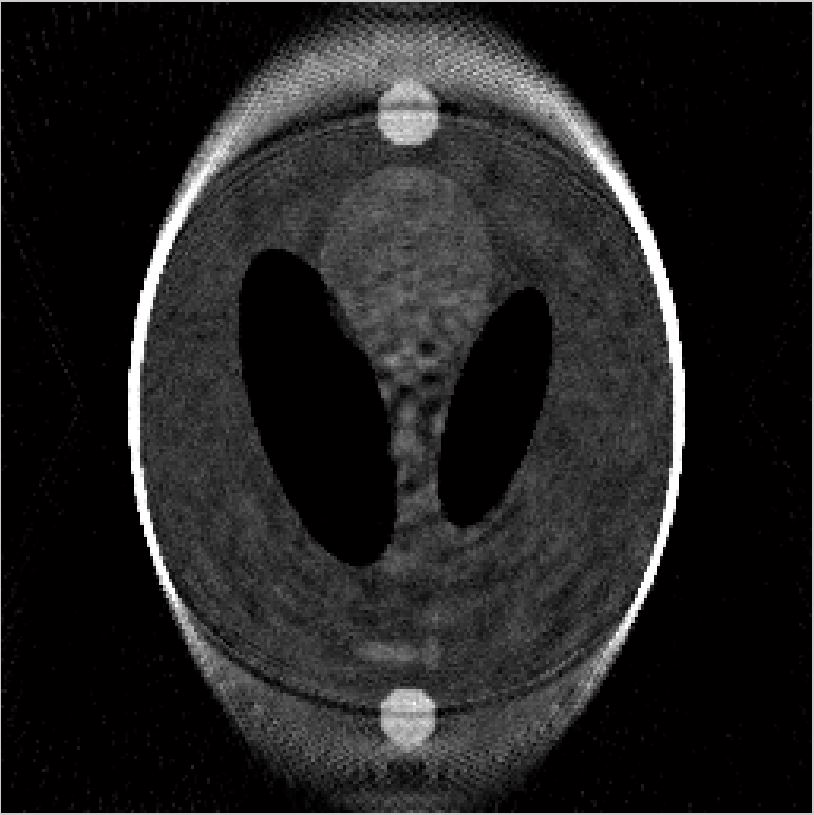}
\end{minipage}&
\begin{minipage}{3cm}
\includegraphics[width=3cm]{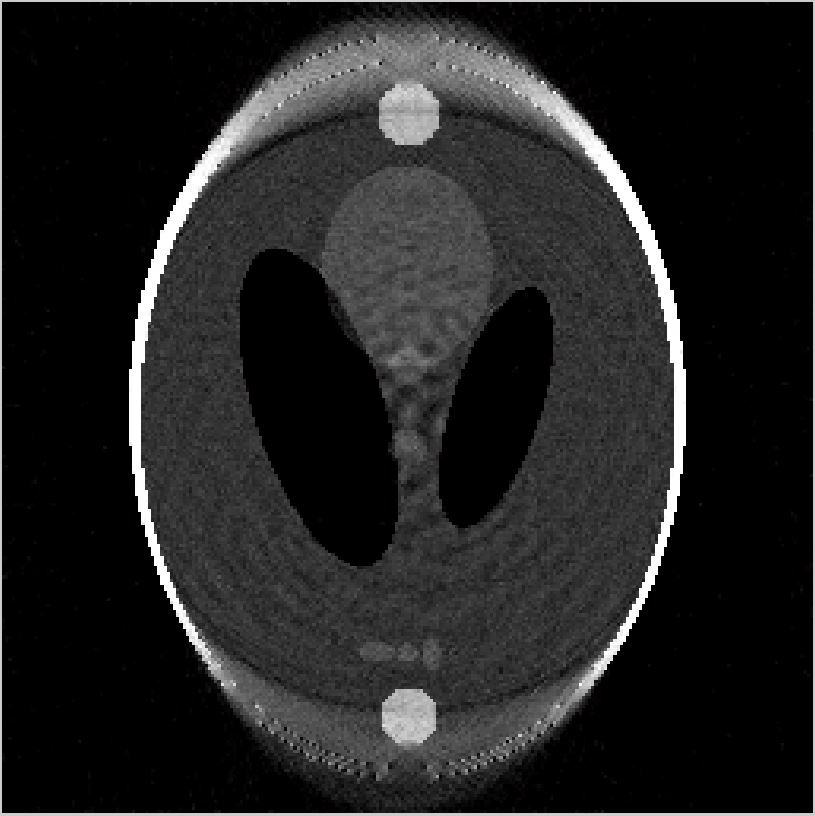}
\end{minipage}&
\begin{minipage}{3cm}
\includegraphics[width=3cm]{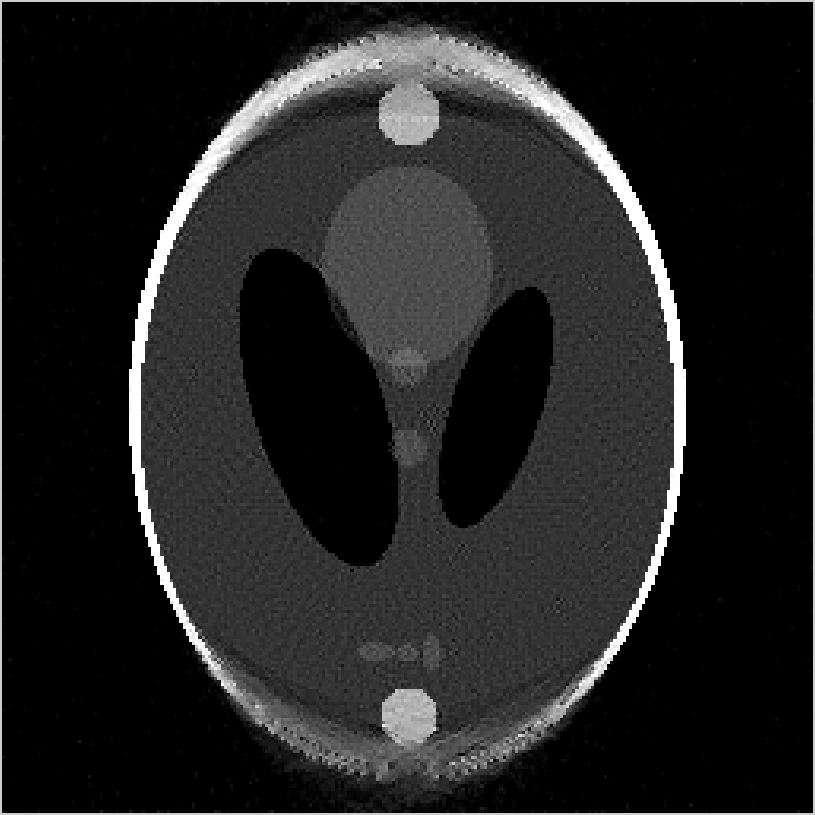}
\end{minipage}
\end{tabular}
\caption{Numerical simulation results.  The top row image is the true image $\wt{\bu}$, followed by rows representing the results using $180$ and $90$ projections respectively.  Images from left to right depicts the reconstructed images obtained by FBP \eref{FBP}, sparsity model based method \eref{Analysis_LT}, wavelet frame based model, and data driven frame based model, respectively.  In this figure, all images are shown in the window level $[0,1]$ for the fair comparison.}\label{Fig41}
\end{center}
\end{figure}

In this section, we compare the wavelet frame based model, the data driven frame based model, and the sparsity model based method \eref{Analysis_LT} with $\W$ being the linear B-spline framelet with $1$ level decomposition.  For \eref{DTDVLT}, we set $\lambda_1=100$ and $\lambda_2=0.01$ for both wavelet frame based model and data driven frame based model.  We test these models using $256\times256$ Shepp-Logan phantom image modified by adding two round shaped objects as the true image $\wt{\bu}$ with the window level $[0,1]$, as shown in the top row of Fig. \ref{Fig41}.  Then we generate the limited projection data $\bsf_0$ through
\begin{align*}
\bsf_0=\bR_{\La}\bsP\wt{\bu}+\eep
\end{align*}
where $\La$ is the discretization of $[0,\pi)\times[-1/2,1/2]$ and $\eep$ is the additive Gaussian noise with zero mean and standard deviation $\s$ being $0.1\%$ of $\max|\bsP\wt{\bu}|$.  The size of discrete Radon transform matrix $\bsP$ depends on both the size of $\wt{\bu}$ and the number of projections $N_p$.  In our simulation, we used $180$ equally spaced projections and $90$ equally spaced projections, both of which are equispatially sampled from $0$ to $\pi$ as two different cases.

Figure \ref{Fig41} describes the reconstructed images for all the tested models and cases.  Compared with FBP, both \eref{Analysis_LT} and \eref{DTDVLT} can recover the image which is missing due to the truncation of the projection data.  However, we can observe that the image reconstructed by our proposed model \eref{DTDVLT} has much clearer geometry information compared with the sparsity model based method \eref{Analysis_LT}.  Besides, compared with the wavelet frame based method, the data driven frame based method recovers the image outside the ROI region.

To provide a quantitative comparison on the models for each tested case, we calculate the PSNR defined as
\begin{align*}
\mathrm{PSNR}:=-20\log_{10}\f{\|\bu-\wt{\bu}\|_2}{N}
\end{align*}
where $N$ denotes the total number of pixels.  In addition, since the goal of our model is to recover the image outside ROI part which is missing due to the limited knowledge of projection, we also calculate the mean structural similarity (MSSIM) defined in \cite{Z.Wang2004}, one of the methods widely used to measure the similarity between the two given images.

\begin{table}[h]
\begin{center}
\begin{tabular}{|c|c|c|c|c|}
\hline
$N_p$ & FBP & Sparsity Model & Wavelet Frame & Data Driven Frame\\
\hline
$180$ & $13.9467$ & $20.848$ & $23.9691$ & $24.1207$\\
\hline
$90$ & $13.8511$ & $19.038$ & $20.8845$ & $22.4786$\\
\hline
\end{tabular}
\caption{Comparison of PSNR of the algorithms.  In both $N_p=180$ and $N_p=90$ cases, the data driven frame based model yields the reconstructed image with highest PSNR, followed by the wavelet frame based model and the sparsity model based method.}\label{Table41}
\end{center}
\end{table}

\begin{table}[h]
\begin{center}
\begin{tabular}{|c|c|c|c|c|}
\hline
$N_p$ & FBP & Sparsity Model & Wavelet Frame & Data Driven Frame\\
\hline
$180$ & $0.4067$ & $0.6776$ & $0.7448$ & $0.7757$\\
\hline
$90$ & $0.2407$ & $0.6466$ & $0.7575$ & $0.8033$ \\
\hline
\end{tabular}
\caption{Comparison of MSSIM of the algorithms.  In both $N_p=180$ and $N_p=90$ cases, the data driven frame based model yields the reconstructed image with highest PSNR, followed by the wavelet frame based model and the sparsity model based method, as in the PSNR.}\label{Table42}
\end{center}
\end{table}

Table \ref{Table41} and \ref{Table42} respectively depict the PSNR and the MSSIM values of the tested models for each cases.  At first, we can see that compared with the FBP, all reconstructed methods yield the reconstructed images with better quantities both in PSNR and MSSIM.  Nevertheless, it can be shown from the tables that our proposed method \eref{DTDVLT} outperforms the sparsity model based method \eref{Analysis_LT}.  Moreover, we can see that the image reconstructed by the data driven frame based method is the closest to the original image in the sense of MSSIM, because at each step we learn tight frames $\W_1$ and $\W_2$ that best represent $\bsf$ and $\bu$ respectively.

Finally, we also observe that the data driven frame based model can recover the image if the number of projections is reduced to $90$.  In the previous researches on the limited tomography, it is in general to increase the number of projections \cite{L.C.B.dosSantos2014}.  However, it is mentioned in \cite{L.C.B.dosSantos2014} that the ambiguity and/or the artifact due to the truncated projection data still remains in the reconstructed image.  Moreover, from the clinical viewpoint, it is not desirable to increase the number of projections too much because of the X-ray radiation dosage.  On the other hand, our proposed data driven frame based method can work very well with the smaller number of projections.  This implies that our method can be helpful for the reduction of X-ray radiation dosage, one of the most challenging issues in the current field of X-ray CT.

However, our method has a limitation.  During the numerical simulation, it has been observed that even though large $\lambda_1$ and $\lambda_2$ can make the reconstructed image more noise-free, but it fails to recover the missing parts of the CT image.  On the other hand, small $\lambda_1$ and $\lambda_2$ can make the reconstructed image noisy, but it performs much better in recovering the image outside of the ROI.  In addition, it requires a mathematical analysis to which extent of data truncation we can reconstruct the entire CT image uniquely and stably.  Even though our reconstruction method via tight frame and simultaneous sinogram extrapolation method has shown a promising result of recovering the missing part of image due to the truncated projection data, further researches will be needed to improve the quality of reconstructed images.

\section{Conclusions}

This paper presents a method of image reconstruction via tight frame and simulataneous sinogram extrapolation for limited tomography.  The numerical simulation shows that our proposed model \eref{DTDVLT} with tight frames presents a possibility of reconstructing the entire image from the limited knowledge of the projection data.  In addition, the missing part of an image due to the truncated projection data can be recovered if we reduce the number of projection angles.  Future direction of the research will focus on the improvement of the reconstructed CT image quality as well as the theoretical condition on the reconstruction of entire image using the limited sinogram.

\section*{Appendix}
\appendix

Provided here is the brief introduction on the concept of tight wavelet frames and data driven tight frames.  For more details, one may refer to \cite{Daubechies1992,I.Daubechies2003,B.Dong2013,Shen2010} for tight wavelet frames and \cite{J.F.Cai2014} for the data driven tight frames.

\subsection*{Tight Wavelet Frames}

A countable set $\msX\subseteq L^2(\R)$ is called a tight frame of $L^2(\R)$ if it satisfies
\begin{align*}
\|f\|_{L^2}^2=\sum_{\varphi\in\msX}\left|\la f,\varphi\ra\right|^2~~~~~~~\text{for all}~~~f\in L^2(\R)
\end{align*}
where $\la\cdot,\cdot\ra$ denotes the inner product in $L^2(\R)$.  Equivalently, every $f\in L^2(\R)$ can be expressed as the linear combination of $\varphi\in\msX$:
\begin{align*}
f=\sum_{\varphi\in\msX}\la f,\varphi\ra\varphi,
\end{align*}
with $\la f,\varphi\ra$ being called the canonical coefficient of $f$.

$\msX$ is called a tight wavelet frame if there exists a finite set $\Ps=\{\psi_1,\cdots,\psi_r\}$ such that
\begin{align*}
\msX=\msX(\Ps)=\{\psi_{l,j,k}(x)=2^{j/2}\psi_l(2^jx-k):1\leq l\leq r~\&~j,k\in\Z\},
\end{align*}
where $\psi_1,\cdots,\psi_r$ are called the framelets.  The compactly supported framelets can be generated by applying the unitary extension principle \cite{A.Ron1997}.  For $n\geq2$, the tight wavelet frame system of $L^2(\R^n)$ can be easily constructed by the tensor products of one dimensional framelets \cite{J.F.Cai2015,J.F.Cai2012,Daubechies1992,B.Dong2013}.

In a discrete setting, a two dimensional image $\bu$ can be regarded as a two dimensional array \cite{J.F.Cai2015,J.F.Cai2012,B.Dong2013a}.  Let $\W$ denote a discrete framelet decomposition (or analysis operator) and let $\W^T$ be the fast reconstruction (or synthesis operator).  Since two dimensional framelets are generated by the tensor products of univariate framelets satisfying the unitary extension principle \cite{A.Ron1997}, it follows that
\begin{align*}
\bu=\W^T\W\bu~~~~~\text{for any image}~~\bu.
\end{align*}
Finally, an $L$-level framelet decomposition of $\bu$ will be denoted as
\begin{align}
\W\bu=\{\W_{l,\bsi}\bu:0\leq l\leq L-1,~0\leq i_1,i_2\leq r\}
\end{align}
where $\bsi=(i_1,i_2)$ is the index of all framelet bands.  If $L=1$, we write $\W_{\bsi}\bu$ for $\W_{l,\bsi}\bu$.

\subsection*{Data Driven Tight Frames}

Even though tensor product framelet based approaches are simple to implement and able to achieve sparse representation of a given image, their major disadvantage lies in that these framelets mostly focus on the horizontal and vertical discontinuities \cite{J.F.Cai2014,J.Liang2014}.  For example, when the discontinuities of an image have complex geometries, the tensor product framelet representation may not be sparse enough \cite{J.F.Cai2014}.

Based on the concept of adaptivity explored by machine learning approaches \cite{M.S.Lewicki2000,K.Kreutz-Delgado2003,M.Elad2006,J.Mairal2008,J.Mairal2009,W.Dong2011}, the goal of data driven tight frame is to construct an adaptive discrete tight frame $\W$ generated by real valued filters $\{\a_1,\cdots,\a_m\}$ for a given image $\bu$.  This is achieved by solving the following minimization problem \cite{J.F.Cai2014}:
\begin{align}\label{DTDVTF}
\min_{\bsd,\{\a_j\}_{j=1}^m}~\|\bsd-\W(\a_1,\cdots,\a_m)\bu\|_2^2+\lambda_1^2\|\bsd\|_0~~~~~\text{s.t.}~~\W^T\W=\I
\end{align}
where $\|\bsd\|_0$ is the number of nonzero elements in the coefficient vector $\bsd$, $\W(\a_1,\cdots,\a_m)$ is the analysis operator generated by filters $\a_1,\cdots,\a_m$, and $\I$ is the identity operator.

The minimization problem \eref{DTDVTF} can be solved by updating $\bsd$ and $\{\a_j\}_{j=1}^m$ alternatively \cite{J.F.Cai2014}.  More precisely, let $\{\a_j^0\}_{j=1}^m$ be the initial filters such as the B-spline framelet filters in \cite{J.F.Cai2015,J.F.Cai2012,B.Dong2013a}.  Then for $k=0,1,\cdots,K-1$, we have
\begin{enumerate}
\item Given filters $\{\a_j^k\}_{j=1}^m$, we define $\W^k=\W(\a_1^k,\cdots,\a_m^k)$.  Then we solve
\begin{align}\label{d-subproblem}
\bsd^k=\arg\min_{\bsd}~\|\bsd-\W^k\bu\|_2^2+\lambda^2\|\bsd\|_0.
\end{align}
\item Given $\bsd^k$, update $\{\a_j^{k+1}\}_{j=1}^m$ by solving the following constrained minimization problem:
\begin{align}\label{filter-update}
\{\a_j^{k+1}\}_{j=1}^m=\arg\min_{\{\a_j\}_{j=1}^m}~\|\bsd^k-\W(\a_1,\cdots,\a_m)\bu\|_2^2~~~~\text{s.t.}~~\W^T\W=\I.
\end{align}
\end{enumerate}
Note that \eref{d-subproblem} has a closed form solution given by the hard threshold:
\begin{align*}
\bsd^k=\bsT_{\lambda}(\W^k\bu)~~~~\text{where}~~~\bsT_{\lambda}(\bsv)[\bk]=\left\{\begin{array}{ccl}
\bsv[\bk]~&\text{if}&~|\bsv[\bk]|\geq\lambda\vspace{0.2em}\\
0~&\text{if}&~|\bsv[\bk]|<\lambda.
\end{array}\right.
\end{align*}
For \eref{filter-update}, the authors in \cite{J.F.Cai2014} has shown that under a proper assumption, the optimal solution has a closed form solution, which simplifies the computation significantly \cite{J.Liang2014}.

\medskip
\noindent{\bf Acknowledgements.}  This work was  partially supported by NSFC (No. 91330102 and GZ1025) and 973 program (No. 2015CB856004).

\bibliographystyle{plain}

\end{document}